\documentclass[preprints,sensors,article,submit,moreauthors,pdftex,12pt,a4paper]{Definitions/mdpi}
\firstpage{1}
\pubvolume{1}
\issuenum{1}
\articlenumber{0}
\pubyear{2021}
\copyrightyear{2020}
\datereceived{}
\dateaccepted{}
\datepublished{}
\hreflink{https://dx.doi.org/} 
\pdfoutput=1
\renewcommand{\linenumbers}{}
\usepackage{amssymb}
\usepackage{enumitem}
\usepackage{mdwtab}
\usepackage{array}
\usepackage{mdwtab}
\usepackage{balance}
\usepackage{bm}

\usetikzlibrary{arrows}

\newcommand{\CommaBin}{\mathbin{\raisebox{0.5ex}{,}}}

\Title{VELOCITY-FREE ATTITUDE CONTROL OF QUADROTORS: \\
A NONLINEAR NEGATIVE IMAGINARY APPROACH}

\TitleCitation{Velocity-Free Attitude Control of Quadrotors: \\
A Nonlinear Negative Imaginary Approach}

\Author{Ahmed~G.~Ghallab $^{1,*}$, and Ian~R.~Petersen $^{1}$}

\AuthorNames{Ahmed~G.~Ghallab, and Ian~R.~Petersen}

\AuthorCitation{Ghallab, A.; Petersen, I.}

\address{%
$^{1}$ \quad Research School of Engineering, The Australian National University, Canberra, ACT 2601, Australia; ahmed.ibrahim@anu.edu.au, i.r.petersen@gmail.com}

\corres{Correspondence: ahmed.ibrahim@anu.edu.au}

\abstract{In this paper, we propose a new approach to the attitude control of quadrotors, by which angular velocity measurements or a model-based observer reconstructing the angular velocity are not needed. The proposed approach is based on recent stability results obtained for nonlinear negative imaginary systems. In specific,  by constructing the respective storage functions, we establish the nonlinear negative imaginary properties of the whole quadrotor system and the quadrotor rotational subsystem. Then, an inner-outer loop method will be implemented to synthesize a strictly negative imaginary controller. This guarantees the robust asymptotic stability of the attitude of the quadrotor about its reference signal in the face of modeling uncertainties and external disturbances.}

\keyword{Nonlinear negative imaginary systems, Quadrotors, Attitude control, Robustness, Feedback stability.)}

\begin{document}
\section{Introduction}
In recent decades, unmanned aerial vehicles (UAVs) have seen increasing interest within the research communities and industry due to their potential for numerous applications including, for instance, inspection, surveillance, data acquisition, and military applications. Their potential future applications include search and rescue, border patrol, surveillance of wildfires, surveillance of traffic and land surveys. An important type of UAVs are quadrotors, which have useful properties such as a simple structure, and low operation and manufacturing costs \cite{du2017distributed, ozbek2016feedback}. In the literature, different control methods have been implemented in order to tackle the quadrotor stability problem such as backstepping \cite{bouadi2007modelling}, sliding mode control \cite{li2014adaptive},  distributed artificial neural networks techniques \cite{tran2020distributed}, and model predictive control \cite{alexis2012model}.

In the existing literature, a great deal of research studies are focused on the attitude control of quadrotor systems which is deemed relatively simple to control. In fact, designing an attitude controller yielding reliability and high level of performance is crucial for many of real-world applications of quadrotors. In most of these attitude control studies, angular velocity measurement is required. However, in practice, velocity measurements in a quadrotor can be noisy or even not available \cite{lizarralde1996attitude}. Alternatively,  observer-based methods have been proposed to reconstruct the angular velocity which lead to inaccurate estimation of the velocity. The inaccurate quality of the velocity measurements can significantly degrade the attitude control performance of such methods.





A widely used approach to tackle the problem of quadrotor attitude control uses the passivity structure of the rigid-body combined with a quaternion-based technique to asymptotically stabilize the attitude of the rigid-body  \cite{lizarralde1996attitude}. However, using quaternions in modeling quadrotors is conceptually challenging and not very intuitive. Motivated by the forgoing, we aim in this paper to apply the nonlinear negative imaginary-based approach, which is recently developed in \cite{ghallab2018extending}  to deal with the robust stability of quadrotor systems that are passive from the input to the derivative of the output (rather than the output as in the classical passivity theory). The nonlinear negative imaginary approach will enable the design of a velocity-free attitude controller by direct use of the Euler angles instead of the quaternion representation of these angles as in \cite{lizarralde1996attitude}.

The rest of the paper is organized as follows: In \textbf{Section 2}, we review some stability results from the NI literature. In \textbf{ Section 3}, we use the Euler-Lagrange dynamics of quadrotors systems to establish the nonlinear negative imaginary property of these systems. Finally in \textbf{Section 4 }, we use an inner-outer loop technique to design a velocity-free attitude controller for the quadrotor system.  We end with some concluding remarks.

\section{Preliminaries}
Negative imaginary (NI) systems theory has been introduced in  \cite{petersen2010} for the control of flexible structures with colocated force actuators and position sensors. In this section, we review some of the related definitions and results from the NI literature in both the linear and nonlinear cases.
\subsection{Negative imaginary systems: Linear Case}
We consider here the following linear time invariant (LTI) system:
\begin{align}
\label{eq:xdot1}
&\dot{x}(t) = Ax(t)+Bu(t), \\
\label{eq:y1} &y(t) = Cx(t)+Du(t)
\end{align}%
where the matrices $A \in \mathbb{R}^{n \times n}, B \in \mathbb{R}^{n \times m}$,  $C\in \mathbb{R}^{m \times n}$, and $D\in \mathbb{R}^{m \times m}$. Assume that the system \eqref{eq:xdot1} and \eqref{eq:y1} has the $m\times m$ real-rational proper  transfer function $G(s):=C(sI-A)^{-1}B+D$.
The frequency domain characterization of the NI property of the LTI system \eqref{eq:xdot1} and \eqref{eq:y1} is given in the following definition.
\begin{Definition} \cite{lanzon2008}. 
A square transfer function matrix $G(s)$ is called negative imaginary if the following conditions are satisfied:
\begin{enumerate}
\item $G(s )$ has no pole at the origin and in $\Re[s]>0$;
\item For all $\omega >0$, such that $j\omega$ is not a pole of $G(s )$, and $j\left( G(j\omega )-G(j\omega )^{T }\right) \geq 0$;
\item If $j\omega_{0}$; $\omega_0\in(0,\infty)$, is a pole of $G(j\omega )$, it is at most a simple pole and the residue matrix $K_{0}= \lim_{ s\rightarrow j\omega_{0}}(s-j\omega_{0})sG(s)$ is positive semidefinite Hermitian.
\end{enumerate}
\end{Definition}

A linear time invariant system of the form \eqref{eq:xdot1} and \eqref{eq:y1} is NI if its transfer function is NI. The definition of the strict negative imaginary property is given as follows.

\begin{Definition} \label{SNI_def}\cite{lanzon2008}. A square transfer function matrix $G(s)$ is strictly negative imaginary if:
\begin{enumerate}
\item  $G(s)$ has no poles in $\Re[s]\geq0$;
\item  $j[G(j\omega)-G^{T}(j\omega)]>0$ for $\omega\in(0,\infty)$.
\end{enumerate}
\end{Definition}
The state-space characterization of strictly negative imaginary systems is given in the following lemma.
\begin{Lemma}\label{SNI}\cite{xiong21010jor}.
Let $(A,B,C,D)$ be a minimal state-space realization of the transfer function matrix $G(s)$. Then $G(s)$ is strictly negative imaginary if and only if:
\begin{enumerate}
\item  $\det(A)\neq 0$, $D=D^{T}$;
\item  there exists a matrix $P=P^{T}>0, P \in \mathbb{R}^{n \times n}$, such that
$$
AP^{-1}+P^{-1}A^{T}\leq 0, \quad \text{and } \quad B+AP^{-1}C^{T}=0;
$$
\item the transfer function matrix $M(s)\backsim
\begin{bmatrix}
\begin{array}{c|c}
A & B \\ \hline LPA^{-1} & 0
\end{array}
\end{bmatrix}$
 has full column rank at $s=jw$ for any $\omega\in(0,\infty)$ where $L^{T}L=-AP^{-1}-P^{-1}A^{T}$. That is, rank $M(j\omega)=m$ for any $\omega\in(0,\infty)$.
\end{enumerate}
\end{Lemma}

\subsection{Nonlinear Negative Imaginary Systems}
We consider the following multi-input multi-output (MIMO) general nonlinear system of the form
\begin{align}
\label{x}\dot{\mathbf{x}}&= \mathbf{f}(\mathbf{x},\mathbf{u})\\
\label{y} \mathbf{y}&= \mathbf{h}(\mathbf{x})
\end{align}
where $\boldsymbol{f}:\mathbb{R}^n\times \mathbb{R}^m\rightarrow\mathbb{R}^n$ is Lipschitz continuous function and $\mathbf{h}:\mathbb{R}^n\rightarrow\mathbb{R}^m$ is continuously differentiable function. The following  definitions give a time-domain characterization of the Negative Imaginary property and its strict version for the class of nonlinear systems of the form \eqref{x}, and \eqref{y}.
\begin{Definition} \cite{ghallab2018extending}\label{nonlinear_NI} Consider the nonlinear system \eqref{x}-\eqref{y}. The system \eqref{x}, \eqref{y} is said to be nonlinear Negative Imaginary if there exists a non-negative function $V:\mathbb{R}^n\rightarrow\mathbb{R}$ of a class $C^1$ such that the following dissipative inequality
\begin{equation}\label{lyapD}
    \dot{V}(\mathbf{x}(t))\leq \dot{\mathbf{y}}^T(t)\mathbf{u}(t),
\end{equation}
holds for all $t \geq 0$. Here, the function is called the storage function.
\end{Definition}
\begin{Definition}\cite{ghallab2018extending}
The system \eqref{x}, \eqref{y} is said to be  a marginally strictly nonlinear NI system if the dissipative inequality \eqref{lyapD} is satisfied, and for all $\mathbf{u}$ and $\mathbf{x}$ such that
\begin{equation}\label{06}
   \dot{V}(\mathbf{x})=\dot{\mathbf{y}}^T(t)\mathbf{u}(t)
\end{equation}
for all $t>0$, then $\lim_{t\rightarrow\infty}\mathbf{u}(t)=\mathbf{0}$.
\end{Definition}
The next definition is a nonlinear generalization of the definition of strictly negative imaginary systems.
\begin{Definition}\cite{ghallab2018extending}
  The system \eqref{x}, \eqref{y} is said to be weakly strictly nonlinear NI system if it is marginally strictly nonlinear NI and globally asymptotically stable with $\mathbf{u}=\mathbf{0}$.
\end{Definition}

In the following, we highlight the main robust stability result introduced in \cite{ghallab2018extending} for the positive feedback interconnection of two nonlinear NI systems. This nonlinear stability result will be used later to stabilize the attitude of a quadrotor system.

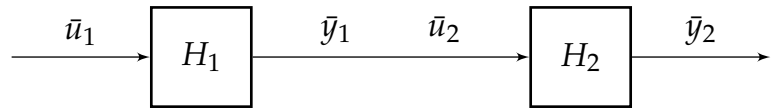
\begin{figure}[H]
\centering
\captionsetup{justification=centering}
\tikzstyle{block} = [draw, thick, rectangle,
    minimum height=2em, minimum width=4em]
\tikzstyle{sum} = [draw, circle,inner sep=0pt,minimum size=1pt, node distance=1cm]
\tikzstyle{input} = [coordinate]
\tikzstyle{output} = [coordinate]
\tikzstyle{pinstyle} = [pin edge={to-,thin,black}]
\tikzstyle{int}=[draw, thick, rectangle, minimum height = 3em,
    minimum width = 3em]
\resizebox{10cm}{!}{%
\begin{tikzpicture}[node distance=2.5cm,auto,>=latex']
\tikzstyle{block} = [draw, thick, rectangle,
    minimum height=3em, minimum width=6em]
    \node [int] (a) {$H_1$};
    \node (b) [left of=a,node distance=2cm, coordinate] {a};
    \node [int] (c) [right of=a,  node distance=4cm] {$H_2$};
    \node [coordinate] (end) [right of=c, node distance=2cm]{};
    \path[->] (b) edge node {$\bar{u}_1$} (a);
    \path[->] (a) edge node {$\bar{y}_1  \quad \quad  \ \bar{u}_2$} (c);
    \draw[->] (c) edge node {$\bar{y}_2$} (end) ;
\end{tikzpicture}}
\caption{Open-loop interconnection of $H_1$ and $H_2$.\label{sys13}}
\end{figure}

Now, consider the following general two MIMO nonlinear systems described by:
\begin{align*}\label{sys}
  H_1{:} \quad \dot{\mathbf{x}}_1 & = \mathbf{f}_1({\mathbf{x}}_1, \mathbf{u}_1)\\
                                \mathbf{y}_1 & = \mathbf{h}_1({\mathbf{x}}_1)  \\
                        \intertext{and}
 H_2{:} \quad  \dot{\mathbf{x}}_2 & =  \mathbf{f}_2({\mathbf{x}}_2, \mathbf{u}_2)\\
                                \mathbf{y}_2 & = \mathbf{h}_2({\mathbf{x}}_2)\\
\end{align*}
where $\mathbf{h}_i:\mathbb{R}^n\rightarrow\mathbb{R}^n$ is a ${C}^1$ function with $\mathbf{h}_i(\mathbf{0})={0}$, $\mathbf{f}_i:\mathbb{R}^n\times \mathbb{R}^n\rightarrow\mathbb{R}^n$ is continuous and locally Lipschitz in ${\mathbf{x}}_i$ for bounded $\mathbf{u}_i$, and where $\mathbf{f}_i(\mathbf{0},\mathbf{0})=\mathbf{0}$. We shall consider the open-loop interconnection of the systems $H_1$ and $H_1$ as shown in Figure~\ref{sys13}. This interconnected system determines the stability properties of the closed-loop system, see \cite{ghallab2018extending}.

We have the following assumptions for the open-loop interconnection of systems $H_1$ and $H_2$.

\begin{Assumption}\label{I}
For any constant $\mathbf{\bar{u}_1}$, there exists a unique solution $\mathbf{\bar{x}_1}, \mathbf{\bar{y}_1}$ to the equations
 \begin{align*}
   \mathbf{0}  &= \mathbf{f}_1(\bar{\mathbf{x}}_1,\bar{\mathbf{u}}_1)\\
   \bar{\mathbf{y}}_1&= \mathbf{h}_1(\bar{\mathbf{x}}_1)
 \end{align*}
 such that $\bar{\mathbf{u}}_1\neq \mathbf{0}$ implies $\bar{\mathbf{x}}_1\neq \mathbf{0}$ and the mapping $\bar{\mathbf{u}}_1 \mapsto \bar{\mathbf{x}}_1$ is continuous.
\end{Assumption}
\vspace{.25cm}
\begin{Assumption}\label{II}
For any constant $\bar{\mathbf{{u}}}_2$, there exists a unique solution $(\bar{\mathbf{x}}_2, \bar{\mathbf{y}}_2)$ to the equations
\begin{align*}
 \mathbf{0}  &= \mathbf{f}_2(\bar{\mathbf{x}}_2,\bar{\mathbf{{u}}}_2)\\
 \bar{\mathbf{y}}_2 &= \mathbf{h}_2(\bar{\mathbf{x}}_2)
\end{align*}
such that $\bar{\mathbf{{u}}}_2\neq \mathbf{0}$ implies $\bar{{\mathbf{x}}}_2\neq \mathbf{0}$.
\end{Assumption}
\vspace{.25cm}
\begin{Assumption}\label{III}
   $\mathbf{h}_1^T(\bar{{\mathbf{x}}}_{1})\mathbf{h}_2(\bar{{\mathbf{x}}}_{2})\geq 0$, for any constant $\bar{\mathbf{{u}}}_1$ where $\bar{\mathbf{{u}}}_2=\bar{\mathbf{y}}_1$.
    \vspace{.25cm}
\end{Assumption}
\begin{Assumption}\label{IV}
    For any constant $\bar{\mathbf{{u}}}_1$, let $(\bar{{\mathbf{x}}}_1, \bar{\mathbf{y}}_1)$ be defined as in Assumption \ref{I} and $(\bar{{\mathbf{x}}}_2, \bar{\mathbf{y}}_2)$ be defined as in Assumption \ref{II} where $\bar{\mathbf{u}}_2=\bar{\mathbf{y}}_1$. Then there exits a constant $0<\gamma<1$ such that for any $\bar{\mathbf{{u}}}_1$ and with $\bar{\mathbf{y}}_2$ defined as in Assumption \ref{II}, the following sector bound condition:
  \begin{equation}\label{sector_bound}
  \bar{\mathbf{y}}_2^T\bar{\mathbf{y}}_2\leq \gamma^2 \bar{\mathbf{{u}}}_1^T\bar{\mathbf{u}}_1,
\end{equation}
holds.
\end{Assumption}

The stability robustness of the positive feedback interconnections of systems $H_1$ and $H_2$ has been established in \cite{ghallab2018extending} using the Lyapunov theory and the LaSalle's invariance principal according to the following theorem.
\begin{Theorem} \cite{ghallab2018extending}\label{main}
Consider a positive feedback interconnection of systems $H_1$ and $H_2$ where $\mathbf{u}_1=\mathbf{y}_2$,  $\mathbf{u_2=y_1}$. Suppose that the system $H_1$ is nonlinear NI and zero-state observable, and the system $H_2$ is weak strict nonlinear NI. Moreover, suppose that Assumptions \ref{I}-\ref{IV} are satisfied. Then, the equilibrium point $\left(\mathbf{x}_1,\mathbf{x}_2\right)=(\mathbf{0},\mathbf{0})$ of the closed-loop system of $H_1$ and $H_2$ is asymptotically stable.
\end{Theorem}


\section{Quadrotor System}
\subsection{Kinematics Model}

Two reference frames are used to study the quadrotor system: a reference frame fixed to the earth $\{R\}(O, x, y, z)$, and a body-fixed frame $\{R_B\}\{O_B, x_B, y_B, z_B\}$, where $O_B$ is fixed to the center of mass of the quadrotor. $\{R_B\}$ is related to $\{R\}$ by a position vector $\bm{\xi}=[x \ y \ z]^T$, describing the position of the center of gravity in $\{R_B\}$ relative to $\{R\}$ and by a vector of three independent angles, known as Euler angles and denoted by $\bm{\eta} = [\phi \ \theta \ \psi]^T$, which represent roll, pitch, and yaw angles of the quadrotor.
Vectors in the body reference frame can be transformed into vectors in the
earth reference frame. For example, given a force $\mathbf{F}_B$, expressed using the coordinates of the body frame, the force $\mathbf{F}$ expressed in the coordinates of the earth frame is:
\begin{equation}
\mathbf{F}= \mathbf{R}_{B\rightarrow E}\mathbf{F}_B
\end{equation}
where $\mathbf{R}_{B\rightarrow E}$ is the transformation (rotation) matrix given by
\begin{equation}\label{rotation_matrix}
  \mathbf{R}:=\left[
      \begin{array}{ccc}
        C_{\theta}C_{\psi} & C_{\psi}S_{\theta}S_{\phi}- C_{\phi}S_{\psi}& C_{\phi}C_{\psi}S_{\theta}+S_{\phi}S_{\psi} \\
         C_{\theta}S_{\psi} & S_{\theta}S_{\phi}S_{\psi}+ C_{\phi}C_{\psi} & C_{\phi}S_{\theta}S_{\psi}-C_{\psi}S_{\phi} \\
        -S_{\theta} & C_{\theta}S_{\phi} & C_{\theta}C_{\phi} \\
      \end{array}
    \right].
\end{equation}
Here $S_{(\cdot)}$ and $C_{(\cdot)}$ represent the functions $\sin(\cdot)$ and $\cos(\cdot)$, respectively.

\vspace{.5cm}

In a similar way, given the angular velocity vector $\boldsymbol{\omega}=(p, q, r)^T$, where $p$, $q$ and $r$ represent the instantaneous angular velocities around the $x_B$-axis, $y_B$-axis and $z_B$-axis respectively, it is related to the rate of change of the yaw, pitch and roll angles in the following way
\begin{equation}
  \boldsymbol{\omega}=\mathbb{W}_{\bm{\eta}} \dot{\boldsymbol{\eta}}
\end{equation}
where
\begin{equation}
  \mathbb{W}_{\bm{\eta}}:=\left[
                \begin{array}{ccc}
                  1 & 0 & -S_{\theta} \\
                  0 & C_{\phi} & S_{\phi}C_{\theta} \\
                  0 & -S_{\phi} & C_{\phi}C_{\theta} \\
                \end{array}
              \right].
\end{equation}

\begin{Remark}
With the assumption of the small angle approximation, we obtain an equality (i.e. $ \mathbb{W}_{\bm{\eta}}=\bm{I}$) between the Euler rates $\dot{\bm{\eta}} = (\dot{\phi} , \dot{\theta}, \dot{\psi})^T$ and the angular velocity vector $\bm{\omega}=(p, q, r)^T$.
\end{Remark}

\subsection{Euler-Lagrange Model: Nonlinear NI Structure}
An Euler-Lagrange approach is adopted in order to write
the equations which describe the translational and rotational motion of
the quadrotor. The Euler-Lagrange dynamics of the quadrotor are given by
\begin{equation}\label{euler_lagrange}
\frac{d}{d t}\left(\frac{\partial \mathcal{L}}{\partial \dot{\mathbf{q}}}\right)-\frac{\partial \mathcal{L}}{\partial \mathbf{q}}=\mathbf{F}
\end{equation}
where $\mathbf{q}=[{\bm{\xi}}^T \ \bm{\eta}^T]^T=[x, y, z, \phi, \theta, \psi]^T \in \mathbb{R}^6$ is the generalized coordinates vector for the quadrotor, and $\mathbf{F}=[\mathbf{F}_{\bm{\xi}}^T \ \bm{\tau}^T]^T$, where $\mathbf{F}_{\bm{\xi}}$ is the thrust force and $\bm{\tau}$ is total torque. The Lagrangian $\mathcal{L}(\mathbf{q},\mathbf{\dot{q}})$ of the quadrotor as the difference between the total kinetic energy $T$ and the potential energy $U$; that is
\begin{equation}\label{lagrangian}
  \mathcal{L}=T_{trans}+T_{rot}-U.
\end{equation}
Here, $T_{trans}$ represents the translational kinetic energy, and is given by
$$T_{trans}=\frac{1}{2} m \dot{\bm{\xi}}^{T} \dot{\bm{\xi}},$$
where $m$ denotes the whole mass of the quadrotor. The term $T_{rot}$ represents the rotational kinetic energy, and is given by
$$T_{rot}=\frac{1}{2} \dot{\bm{\eta}}^{T} \mathbf{J} \dot{\bm{\eta}},$$
where $\mathbf{J}$ represent the rotational inertia matrix and is defined in the body frame as
\begin{equation} \mathbf{J}=
  \left[
     \begin{array}{ccc}
        J_{x} & 0 & 0 \\
        0 & J_{y} & 0 \\
        0 & 0 & J_{z}
     \end{array}
   \right],
\end{equation}
which is diagonal because of the quadrotor's symmetry.


\begin{specialtable}[t!]
\centering
\captionsetup{justification=centering}
\caption{Quadrotor Parameters.\label{tab1}}
\tablesize{\normalsize} 
\begin{tabular}{ccc}
\toprule
\textbf{Definition}	& \textbf{Parameter} & \textbf{Unit}\\
\midrule
 Quadrotor mass & $m$ & $\mathrm{kg}$ \\
 Gravitational acceleration & $g$ & $\mathrm{m}/\mathrm{s}^{2}$ \\
 Arm length  & $l$ & $\mathrm{m}$ \\
 Thrust coefficient  & $b$ & $\mathrm{N}\cdot\mathrm{s}^{2}/\mathrm{rad}^{2}$ \\
 Drag coefficient  & $d$ & $\mathrm{N}\cdot\mathrm{s}^{2}/\mathrm{rad}^{2}$ \\
 Roll inertia  x-axis  & $J_{x}$ & $\mathrm{kg}\cdot\mathrm{m}^{2}$ \\
 Pitch inertia  y-axis  & $J_{y}$ & $\mathrm{kg}\cdot\mathrm{m}^{2}$ \\
 Yaw inertia  & $J_{z}$ & $\mathrm{kg}\cdot\mathrm{m}^{2}$ \\
\bottomrule
\end{tabular}
\end{specialtable}

In terms of the generalized coordinates, the Lagrangian \eqref{lagrangian} can be written as
\begin{equation}\label{lagrangian2}
  \mathcal{L}(\mathbf{q},\mathbf{\dot{q}})=\frac{1}{2} m \dot{\bm{\xi}}^{T} \dot{\bm{\xi}}+\frac{1}{2} \dot{\bm{\eta}}^{T} \mathbf{J} \dot{\bm{\eta}}-m g z=\frac{1}{2} \dot{\mathbf{q}}^{T} \mathbf{M} \dot{\mathbf{q}}+\mathbf{B} \mathbf{q}
\end{equation}
where
$$
\mathbf{M}=\left[\begin{array}{cc}
m\mathbf{I}_{3\times3} & \mathbf{0}_{\mathbf{3} \times \mathbf{3}} \\
\mathbf{0}_{\mathbf{3} \times \mathbf{3}} & \mathbf{J}
\end{array}\right], \quad \mathbf{B}=\left[0 \ 0 \ \ m g \ \ 0 \ \ 0 \ \ 0\right]^{T}.
$$
Using the Lagrangian \eqref{lagrangian2}, equation \eqref{euler_lagrange} can be written as follows
\begin{equation}\label{lagrange}
\mathbf{M}(\mathbf{q}(t)) \ddot{\mathbf{q}}(t)+\mathbf{C}(\mathbf{q}(t), \dot{\mathbf{q}}(t)) \dot{\mathbf{q}}(t)+\mathbf{\mathbf{G(q)}}=\mathbf{F(t)}
\end{equation}
where $\mathbf{M(q)}$ is the inertia matrix, which is symmetric and positive definite, $\mathbf{C({q}, \dot{q})}$ is the Coriolis and centrifugal matrix, $\mathbf{G(q)}$ is the gravitational vector, and $\mathbf{F}$ is the input (generalized forces) of the system. The the Coriolis and centrifugal matrix is given by
\begin{equation}\label{Coriolis}
\mathbf{C}(\mathbf{\mathbf{q}}, \dot{\mathbf{q}})=\frac{\mathrm{d}}{\mathrm{d} t}\mathbf{M(q)}
-\frac{1}{2} \frac{\partial}{\partial \mathbf{q}}\left(\dot{\mathbf{q}}^{T} \mathbf{M} \right)
\end{equation}

\begin{Remark}
Note that the representation \eqref{Coriolis} of the matrix $\mathbf{C({q}, \dot{q})}$
is not unique \cite{lewis2003robot}. By defining $\mathbf{C}(\mathbf{\mathbf{q}}, \dot{\mathbf{q}})$ using the Christoffel symbols, $\mathbf{M(q)} - 2\mathbf{C}(\mathbf{\mathbf{q}}, \dot{\mathbf{q}})$ is skew-symmetric \cite{ortega1989adaptive}.
\end{Remark}

Next we show in the following lemma that the quadrotor system \eqref{lagrange} with $\mathbf{F}$ as input  and $\mathbf{q}$ as output is nonlinear negative imaginary.

\begin{Lemma}
Consider the system with input $\mathbf{F}$ and output $\mathbf{q}$. Then the system is nonlinear negative imaginary with a positive definite storage function given by
\begin{equation}
    V(\mathbf{q}, \dot{\mathbf{q}})=\frac{1}{2}\dot{\mathbf{q}}^{T} M(\mathbf{q}) \dot{\mathbf{q}}+\mathbf{U(q)}
\end{equation}
where $\mathbf{U(q)}$ is the potential energy where $\mathbf{G(q)}=\frac{\partial \mathbf{U(q)}}{\partial \mathbf{q}}$.  $\mathbf{U(q)}$ is assumed to have an absolute minimum at $\mathbf{q}=\mathbf{0}$.
\end{Lemma}
\vspace{.5em}
\begin{proof}It can be easily shown that $V$ is positive definite. Taking the first time derivative of the function $V$ we obtain
\begin{equation*}
\begin{aligned}
\frac{d V}{d t}(\mathbf{q}, \dot{\mathbf{q}}) &=\dot{\mathbf{q}}^{T} \mathbf{M}(\mathbf{q}) \ddot{\mathbf{q}}+\frac{1}{2} \dot{\mathbf{q}}^{T} \dot{\mathbf{M}}(\mathbf{q}) \dot{\mathbf{q}}+\mathbf{G(\mathbf{q})}\dot{\mathbf{q}} \\
&=\dot{\mathbf{q}}^{T}\left[-\mathbf{C}(\mathbf{q}, \dot{\mathbf{q}}) \dot{\mathbf{q}}-\mathbf{G(q)}+\mathbf{F}\right]\\&\hspace{3cm} +\frac{1}{2} \dot{\mathbf{q}}^{T} \dot{\mathbf{M}}(\mathbf{q}) \dot{\mathbf{q}}+\mathbf{G(q)}\dot{\mathbf{q}} \\
&=\dot{\mathbf{q}}^{T} \mathbf{F}+\frac{1}{2} \dot{\mathbf{q}}^{T}\left[\dot{\mathbf{M}}(\mathbf{q})-2 \mathbf{C}(\mathbf{q}, \dot{\mathbf{q}})\right] \dot{\mathbf{q}}= \dot{\mathbf{q}}^{T} \mathbf{F},
\end{aligned}
\end{equation*}
which shows that the systems is nonlinear negative system from $\mathbf{F}$ to $\mathbf{q}$.
\begin{flushright}
\qedhere
\end{flushright}
\end{proof}
\begin{Remark}
By revealing the nonlinear negative imaginary structure of the quadrotor model, we can leverage powerfull techniques from negative imaginary systems theory to different quadrotors control problems such as, for instance, trajectory tracking and altitude control problems.
\end{Remark}

\subsection{State-Space Model}
Since the stability results of \cite{ghallab2018extending} deal with the state-space representation of  nonlinear systems, we aim to use the equivalent state-space representation of the quadrotor system.  We see that the Lagrangian contains no cross-terms in the kinetic energy combining $\bm{\xi}$ and $\dot{\bm{\eta}}$, so the Euler-Lagrange equation \eqref{euler_lagrange} partitions into two parts; that is, the translational and rotational components.
The translational equation of the quadrotor is described by the following equation
\begin{equation}
m \ddot{\bm{\xi}}+\left(\begin{array}{c}
0 \\
0 \\
m g
\end{array}\right)=\mathbf{F}_{\xi},
\end{equation}
with ${\mathbf{F}}_\xi$ is the thrust force generated by the four rotors and is given by ${\mathbf{F}}_\xi=\sum_{i=1}^4 b \omega_i^2$, where $\omega_i$ is \textit{i}th-rotor's speed and $b$ is the thrust factor.
The rotational subsystem describing the roll,  pitch and yaw rotations  of the quadrotor is described by the following equation
\begin{equation}\label{euler_lagrange_rot}
\mathbf{J}(\bm{\eta}) \ddot{\bm{\eta}}+\frac{d}{d t}\{\mathbf{J}(\bm{\eta})\} \dot{\bm{\eta}}-\frac{1}{2} \frac{\partial}{\partial \bm{\eta}}(\dot{\bm{\eta}}^{T} \mathbf{J}(\bm{\eta}) \dot{\bm{\eta}})=\bm{\tau}
\end{equation}
or, by appropriate definition of variables,
\begin{equation}
\mathbf{J}(\bm{\eta}) \ddot{\bm{\eta}}+C(\bm{\eta}, \dot{\bm{\eta}}) \dot{\bm{\eta}}=\bm{\tau}
\end{equation}
where  the input vector  $\bm{\tau}=[\tau_\phi \ \tau_\theta \ \tau_\psi]^T \in \mathbb{R}^3$  is the total torque in the pitch, roll, and yaw. In matrix form, the vector $\bm{\tau}$ is defined in terms of the four rotor speeds as follows,
\begin{equation}
\bm{\tau} =\left(\begin{array}{c}
l b\left(\omega_{2}^{2}-\omega_{4}^{2}\right) \\
l b\left(\omega_{1}^{2}-\omega_{3}^{2}\right) \\
d\left(\omega_{1}^{2}+\omega_{3}^{2}-\omega_{2}^{2}-\omega_{4}^{2}\right)
\end{array}\right)
\end{equation}
where $l$ is the arm length, the distance from
the axis of rotation of the rotors to the center of the quadrotor and $d$ is the drag force.

By defining the input of the quadrotor as follows
\begin{equation}
\left(\begin{array}{l}
F_{\xi} \\
\tau_{\phi} \\
\tau_{\theta} \\
\tau_{\psi}
\end{array}\right)=\left(\begin{array}{l}
u_{1} \\
u_{2} \\
u_{3} \\
u_{4}
\end{array}\right)=\left(\begin{array}{cccc}
b & b & b & b \\
0 & b & 0 & -b \\
b & 0 & -b & 0 \\
d & -d & d & -d
\end{array}\right)\left(\begin{array}{l}
\omega_{1}^{2} \\
\omega_{2}^{2} \\
\omega_{3}^{2} \\
\omega_{4}^{2}
\end{array}\right),
\end{equation}
we obtain the overall quadrotor dynamic model in the following form
\begin{equation}\label{state_space}
\begin{aligned}
\ddot{x} &=-(\cos \phi \sin \theta \cos \psi+\sin \phi \sin \psi) \cdot \frac{u_{1}}{m} \\
\ddot{y} &=-(\cos \phi \sin \theta \sin \psi-\sin \phi \cos \psi) \cdot \frac{u_{1}}{m} \\
\ddot{z} &=g-(\cos \phi \cos \theta) \cdot \frac{u_{1}}{m} \\
\ddot{\phi} &=\dot{\theta} \dot{\psi}\left(\frac{{J}_y-J_{z}}{J_{x}}\right)-\frac{J_{r}}{J_{x}} \dot{\theta} g(\mathbf{u})+\frac{l}{J_{x}} u_{2} \\
\ddot{\theta} &=\dot{\phi} \dot{\psi}\left(\frac{J_{z}-I_{x}}{J_{y}}\right)+\frac{J_{r}}{J_{y}} \dot{\phi} g(\mathbf{u})+\frac{l}{J_{y}} u_{3} \\
\ddot{\psi} &=\dot{\phi} \dot{\theta}\left(\frac{J_{x}-J_{y}}{J_{z}}\right)+\frac{1}{J_{z}} u_{4}
\end{aligned}
\end{equation}
The system \eqref{state_space} can be represented in the form $\dot{\mathbf{x}}=\mathbf{f}(\mathbf{x}, \mathbf{u})$ with the 12-dimensional state vector $\mathbf{x}=\left(x_{1}, \ldots, x_{12}\right)^{T}=\left(\phi, \dot{\phi} , \theta , \dot{\theta} , \psi , \dot{\psi} , z , \dot{z} , x , \dot{x} , y , \dot{y}
\right)^{T}$, and the input vector $\mathbf{u}=\left(u_{1}, u_{2}, u_{3}, u_{4}\right)^{T}$ :
\begin{equation}
\mathbf{f}(\mathbf{x}, \mathbf{u})=\left(\begin{array}{c}
x_{2} \\
x_{4} x_{6} a_{1}-x_{4} a_{2} g(u)+b_{1} u_{2} \\
x_{4} \\
x_{2} x_{6} a_{3}+x_{2} a_{4} g(u)+b_{2} u_{3} \\
x_{6} \\
x_{4} x_{2} a_{5}+b_{3} u_{4}\\
x_8\\
g-\frac{u_{1}}{m} \cos x_{1} \cos x_{3} \\
x_{10} \\
-\frac{u_{1}}{m}\left(\sin x_{1} \sin x_{5}+\cos x_{1} \sin x_{3} \cos x_{5}\right) \\
x_{12} \\
\frac{u_{1}}{m}\left(\sin x_{1} \cos x_{5}-\cos x_{1} \sin x_{3} \sin x_{5}\right)
\end{array}\right)
\end{equation}
where
\begin{equation*}
\begin{array}{l|l}
a_{1}=\left(J_{y}-J_{z}\right) / J_{x} & b_{1}=l / J_{x} \\
a_{2}=J_{r} / J_{x} & b_{2}=l / J_{y} \\
a_{3}=\left(J_{z}-J_{x}\right) / J_{y} & b_{3}=1 / J_{z} \\
a_{4}=J_{r} / J_{y} \\
a_{5}=\left(J_{x}-J_{y}\right) / J_{z} &
\end{array}
\end{equation*}
\section{Attitude Control Design}
As it can be seen from the previous section,  the state- space model of the quadrotor can be divided into subsystems, one of which, the rotational subsystem, describes the dynamics of the attitude (i.e. the angles) and the other describes the translation of the quadrotor. In this section, we are interested in the problem of stabilizing the attitude of the quadrotor by using a nonlinear negative imaginary approach. For this purpose we confine ourselves to the rotational subsystem whose state is a restriction of $\mathbf{x}$ to the last $6$ components representing the roll, pitch and yaw angles and their time derivatives. The rotational subsystem is then described by $\mathbf{\dot{x}}=\mathbf{f}_{\alpha}(\mathbf{x}, \mathbf{u})$:
\begin{equation}\label{rot_subsystem}
\mathbf{\dot{x}}=\left(\begin{array}{c}
x_{2} \\
x_{4} x_{6} a_{1}-x_{4} a_{2} g(u)+b_{1} u_{2} \\
x_{4} \\
x_{2} x_{6} a_{3}+x_{2} a_{4} g(u)+b_{2} u_{3} \\
x_{6} \\
x_{4} x_{2} a_{5}+b_{3} u_{4}
\end{array}\right).
\end{equation}

The use of a multi-loop control architecture has been proposed in recent studies for a variety of quadrotor control problems; see for instance \cite{bolandi2013attitude, thanh2018simple}. In this section, we propose an inner-outer loop architecture to robustly stabilize the attitude, i.e. the angles, around desired reference signal  $\left(\phi_{d}, \theta_{d}, \psi_{d}\right)^{T}=\left(x_{1}^d, x_{3}^d, x_{5}^d\right)^{T}.$


\subsection{Inner Control Loop}


The inner control loop is mainly designed due to the free motion behavior of the quadrotor. By denoting the vector of error attitude $\bm{\Tilde{\eta}}=\left(\bm{\eta}- \bm{\eta_{d}}\right)$, where $\bm{\eta_{d}}$ is the reference signal, we define the following feedback control law
\begin{equation}
\bm{\tau}=-\mathbf{K_{p}} \bm{\Tilde{\eta}}+\bm{v}
\end{equation}
where $\bm{v}=(v_1,v_2,v_3)^T$ denotes the new input torque of the quadrotor, and $\mathbf{K_{p}}=\operatorname{diag}\left( k_{p}^{\phi},  k_{p}^{\theta}, k_{p}^{\psi}\right)$ is a positive diagonal matrix and the diagonal elements are used as tuning parameters. The designed torque $\bm{\tau}$ is then defined as follows,
\begin{equation}\label{control_law}
\left(\begin{array}{l}
u_2\\ u_3\\u_4
\end{array}\right)=\left(\begin{array}{c}
-k_{p}^\phi\left(x_1-x_1^{d}\right)+ v_1\\
-k_{p}^\theta\left(x_2-x_3^d\right)+ v_2\\
-k_{p}^\psi\left(x_3-x_3^d\right)+ v_3
\end{array}\right)
\end{equation}
Using \eqref{rot_subsystem}, \eqref{control_law}, and setting the desired reference signal $\left(x_{1}^d, x_{3}^d, x_{5}^d\right)^{T}=(0,0,0)^T$, the rotational subsystem can be put in the form $\mathbf{\dot{x}}=\mathbf{\tilde{f}}_{\alpha}(\mathbf{x}, \mathbf{u})$:
\begin{equation}\label{reshaped_rot}
\mathbf{\dot{x}}=\left(\begin{array}{c}
x_{2} \\
x_{4} x_{6} a_{1}-x_{4} a_{2} g(u) -b_{1}k_{p}^\phi\left(x_1-x_1^{d}\right)+ b_1 v_1\\
x_{4} \\
x_{2} x_{6} a_{3}+x_{2} a_{4} g(u) -b_{2}k_{p}^\theta\left(x_2-x_3^d\right)+b_2 v_2 \\
x_{6} \\
x_{4} x_{2} a_{5}-b_{3}k_{p}^\psi\left(x_3-x_3^d\right)
+b_3 v_3
\end{array}\right)
\end{equation}

\begin{figure}[H]
\widefigure
\centering
\includegraphics{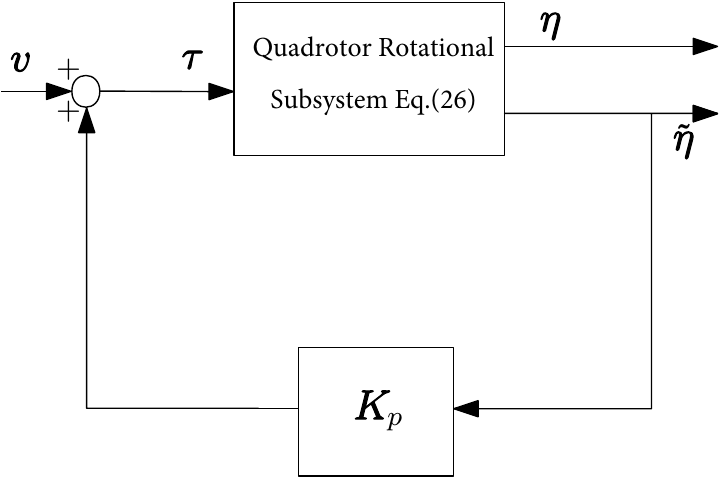}
    \caption{Inner control loop (rotational subsystem \eqref{reshaped_rot}).\label{inner_loop}}
\end{figure}

The above rotational dynamical system, as shown in Figure~\ref{inner_loop}, can be seen as a nonlinear negative imaginary system from the input $\bm{v}$ to the output $\bm{\eta}$ according to the following theorem.

\begin{Theorem}
 Consider the quadrotor rotational subsystem \eqref{reshaped_rot} with  $\bm{v}$ as input and $\bm{\eta}$ as output. Then the system \eqref{reshaped_rot} is a nonlinear negative imaginary system with respect to the following positive-definite storage function,
 \begin{equation}
    V=\frac{1}{2}\dot{\bm{\eta}}^{T} \mathbf{J}(\bm{\eta}) \dot{\bm{\eta}}+\frac{1}{2}\left( \bm{\eta}- \bm{\eta}_d\right)^T \mathbf{K}_p\left(\bm{\eta}-\bm{\eta}_d\right).
\end{equation}
\end{Theorem}
\begin{proof}
It can be easily shown that $V$ is a valid storage function
function; since
the rotational inertia matrix $\mathbf{J}$ is positive definite and $\mathbf{K}_p>0$. Taking the time derivative of $V$ we obtain
\begin{equation*}
\begin{aligned}
\frac{d V}{d t}(\bm{\eta}, \dot{\bm{\eta}}) &=\dot{\bm{\eta}}^{T} \mathbf{J}(\bm{\eta}) \ddot{\bm{\eta}}+\frac{1}{2}\dot{\bm{\eta}}^{T} \mathbf{\dot{J}}(\bm{\eta}) \dot{\bm{\eta}} + \bm{\tilde{\eta}}^T \mathbf{K}_p \bm{\dot{\tilde{\eta}}}\\
&=\dot{\bm{\eta}}^{T}[-\mathbf{C}(\bm{\eta}, \dot{\bm{\eta}}) \dot{\bm{\eta}}+\bm{\tau}]+\frac{1}{2} \dot{\bm{\eta}}^{T} \mathbf{\dot{J}}(\bm{\eta}) \dot{\bm{\eta}} +\bm{\tilde{\eta}}^T \mathbf{K}_p \bm{\dot{\tilde{\eta}}}\\
&=\dot{\bm{\eta}}^{T} \bm{v}+ \frac{1}{2}\dot{\bm{\eta}}^{T}[\mathbf{\dot{J}}(\bm{\eta})-2 \mathbf{C}(\bm{\eta}, \dot{\bm{\eta}})] \dot{\bm{\eta}}= \dot{\bm{\eta}}^{T}\bm{v},
\end{aligned}
\end{equation*}
which shows that the system \eqref{reshaped_rot} is nonlinear negative imaginary from input $\bm{v}$ to output $\bm{\eta}$.
\begin{flushright}
\qedhere
\end{flushright}
\end{proof}
\end{paracol}
\begin{paracol}{2}
\linenumbers
\switchcolumn
\subsection{Outer Control Loop}
We aim here to design a positive feedback closed-loop system that ensure the asymptotic stability of the quadrotor attitude in virtue of Theorem \ref{main}. For simplicity, we will use the following linear MIMO integral resonant controller as the outer control loop controller,
\begin{equation}\label{SNI_controller}
    \mathbf{C}_v(s)=[s\mathbf{I}+\bm{\Gamma}\bm{\Phi}]^{-1}\bm{\Gamma}.
\end{equation}
Here,  $\bm{\Gamma}$ and $\bm{\Phi}$ are positive-definite matrices where  $\bm{\Phi}=\operatorname{diag}\left(\phi, \phi, \phi\right)$, $\phi$ is the tuning parameter. The transfer function matrix $\mathbf{C}_v(s)$ is strictly negative imaginary \cite{petersen2010}. The dc-gain (the gain of the system at steady-state) of the controller is  $\mathbf{C}_v(0)=\bm{\Phi}^{-1}$.

In order to apply Theorem \ref{main}, we need to validate Assumptions \ref{I}-\ref{IV} on the open-loop interconnection (as shown in Figure~\ref{sys3}) of the quadrotor rotational subsystem \eqref{reshaped_rot} and the SNI controller \eqref{SNI_controller}. As it can be seen, in the steady-state case, we have $\mathbf{\tilde{f}}_\alpha(\mathbf{\bar{x}},\mathbf{\bar{u}})=\mathbf{0}$, where for every constant value of $\bm{\bar{v}}$ there is a corresponding unique value of $\mathbf{\bar{x}}$ such that

\begin{equation}\label{steady_state_sol}
\begin{aligned}
 \bar{x}_1&=\frac{\bar{v}_2}{k_{p}^\phi}\CommaBin\\
\bar{x}_3&=\frac{\bar{v}_3}{k_{p}^\theta}\CommaBin\\
\bar{x}_5&=\frac{\bar{v}_5}{k_{p}^\psi}\CommaBin
\end{aligned}
\end{equation}
which shows that Assumption \ref{I} holds. Also, Assumption \ref{II} trivially holds since the controller \eqref{SNI_controller} is an LTI system. We can easily see that Assumption \ref{III} is valid since we have $\mathbf{\bar{y}}_c=\bm{\Phi}^{-1} \mathbf{\bar{u}}_c=\frac{1}{\displaystyle \phi}\mathbf{\bar{y}}$, which yields
$$ \mathbf{\bar{y}}^T\mathbf{\bar{y}}_c=\frac{1}{\phi}\mathbf{\bar{y}}^T\mathbf{\bar{y}}\geq0.$$
\begin{figure}[t!]
\centering
\tikzstyle{block} = [draw, thick, rectangle,
    minimum height=2em, minimum width=4em]
\tikzstyle{sum} = [draw, circle,inner sep=0pt,minimum size=1pt, node distance=1cm]
\tikzstyle{input} = [coordinate]
\tikzstyle{output} = [coordinate]
\tikzstyle{pinstyle} = [pin edge={to-,thin,black}]
\tikzstyle{int}=[draw, thick, rectangle, minimum height = 3em,
    minimum width = 3em]
\resizebox{10cm}{!}{%
\begin{tikzpicture}[node distance=2.5cm,auto,>=latex']
\tikzstyle{block} = [draw, thick, rectangle,
    minimum height=3em, minimum width=6em]
    \node [int] (a) {$\begin{matrix}
   \mathbf{\tilde{f}}_\alpha(\mathbf{\bar{x}},\mathbf{\bar{u}})=0, \\
    \mathbf{\bar{y}}= \bm{\bar{\eta}}
  \end{matrix}$};
    \node (b) [left of=a,node distance=3cm, coordinate] {a};
    \node [int] (c) [right of=a,  node distance=4cm] {$\displaystyle \bm{\Phi}^{-1}$};
    \node [coordinate] (end) [right of=c, node distance=2cm]{};
    \path[->] (b) edge node {$\mathbf{\bar{u}}=\bm{\bar{v}}$} (a);
    \path[->] (a) edge node {$\mathbf{\bar{y}}=\bm{\bar{\eta}}  \quad \quad   \mathbf{\bar{u}}_c$} (c);
    \draw[->] (c) edge node {$\mathbf{\bar{y}}_c$} (end) ;
\end{tikzpicture}}
\caption{Open-loop interconnection (in the steady-state case) of the quadrotor rotational subsystem \eqref{reshaped_rot} and the SNI controller \eqref{SNI_controller} (where'$c$' refers to the controller).\label{sys3}}
\end{figure}
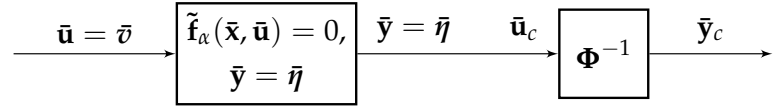
Lastly,  using \eqref{steady_state_sol} we see that
\begin{align*}
\mathbf{\bar{y}}_{c}^{T} \mathbf{\bar{y}}_{c}&=\frac{1}{\phi^{2}}\left[\bar{x}_1+\bar{x}_3+\bar{x}_5\right]\\&=\frac{1}{\phi^{2}}\left[\frac{\bar{v}_{2}^{2}}{\left(\displaystyle k_{p}^{\phi}\right)^{2}} +\frac{\bar{v}_{3}^{2}}{\left(\displaystyle k_{p}^{\theta}\right)^{2}} +\frac{\bar{v}_{4}^{2}}{\left(\displaystyle k_{p}^{\psi}
\right)^{2}} \right]\\
&\leq \frac{1}{\phi^{2}} \max _{i=\phi,\theta,\psi}\left\{\frac{1}{\left(\displaystyle k_{p}^{i}\right)^{2}}\right\} \left(\bar{v}_{2}^{2}+\bar{v}_{3}^{2}+\bar{v}_{4}^{2}\right)\\
&= \frac{1}{\phi^{2}} \max _{i=\phi,\theta,\psi}\left\{\frac{1}{\left(k_{p}^{i}\right)^{2}}\right\} \bm{\bar{v}}^T\bm{\bar{v}}\\
&= \gamma^2\bm{\bar{v}}^T\bm{\bar{v}},
\end{align*}
where $$\gamma^2=\frac{1}{\phi^{2}} \max _{i=\phi,\theta,\psi}\left\{\frac{1}{\left(k_{p}^{i}\right)^{2}}\right\}.$$
To ensure that Assumption \ref{IV} is valid so that $\gamma^2<1$, we choose the outer-loop controller DC-gain should satisfy
\begin{equation}
\phi^{2}>\max _{i=\phi,\theta,\psi}\left\{\frac{1}{\left(\displaystyle k_{p}^{i}
\right)^{2}}\right\}.
\end{equation}

\begin{figure}[t!]
\centering
\includegraphics{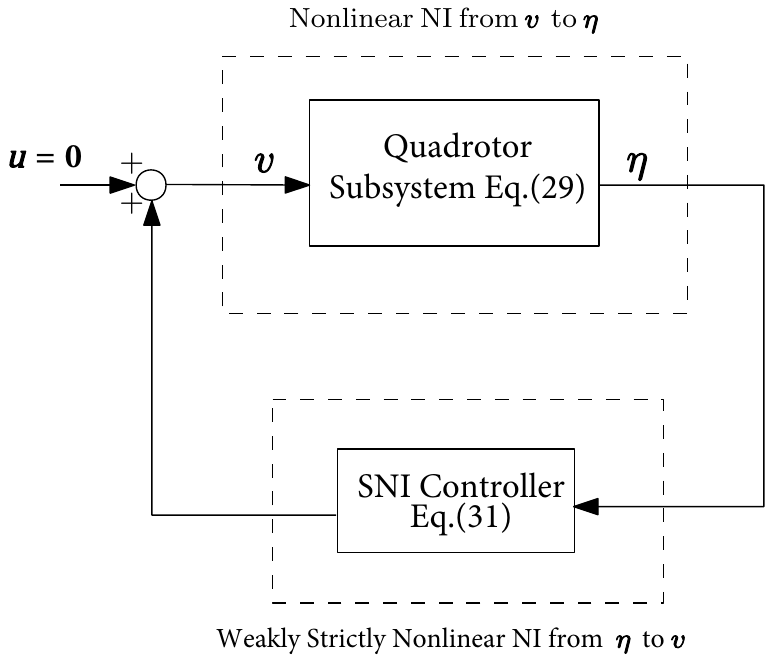}
\caption{Block diagram of the proposed attitude control system comprised of a positive feedback interconnection of the rotational subsystem \eqref{reshaped_rot} and SNI controller \eqref{SNI_controller}.}
\label{fig2}
\end{figure}
Therefore, we have found a lower bound on the DC-gain of the outer controller that is necessary to achieve
asymptotic stability of the attitude vector around the reference signal by virtue of Theorem \ref{main}. We summarize the above result in the following theorem.


\begin{Theorem}
Consider the closed-loop system, as in Figure~\ref{fig2}, of the quadrotor rotational subsystem \eqref{reshaped_rot} and the strictly negative imaginary controller \eqref{SNI_controller}. Assume that the system  \eqref{reshaped_rot} is zero-state observable. Then the closed-loop system is asymptotically stable.
\end{Theorem}
\begin{Remark}
The above stability result leads to a robust control system since the stability is guaranteed irrespective of the quadrotor and controller parameters so long as the condition is satisfied.
\end{Remark}


\section{Conclusions}
The nonlinear negative imaginary systems theory has been used to control the attitude of the quadrotor model around a reference signal. We have proposed a technique based on recently developed results for the stability robustness of positive feedback interconnections of nonlinear NI systems. An inner-outer loop technique has been implemented to design a positive feedback control system that robustly stabilizes the quadrotor's attitude in the face of model uncertainties and disturbances.


\vspace{6pt}




\funding{This work was supported by the Australian Research Council under grant DP160101121.}

\institutionalreview{Not applicable.}

\informedconsent{Not applicable.}

\dataavailability{Not applicable.}


\conflictsofinterest{The authors declare no conflict of interest.}

\end{paracol}


\externalbibliography{yes}
\bibliography{ifacbib2}

\begin{thebibliography}{999}

\bibitem[Du \em{et~al.}(2017)Du, Zhu, Wen, Duan, and L{\"u}]{du2017distributed}
Du, H.; Zhu, W.; Wen, G.; Duan, Z.; L{\"u}, J.
\newblock Distributed formation control of multiple quadrotor aircraft based on
  nonsmooth consensus algorithms.
\newblock {\em IEEE transactions on cybernetics} {\bf 2017}, {\em
  49},~342--353.

\bibitem[{\"O}zbek \em{et~al.}(2016){\"O}zbek, {\"O}nkol, and
  Efe]{ozbek2016feedback}
{\"O}zbek, N.S.; {\"O}nkol, M.; Efe, M.{\"O}.
\newblock Feedback control strategies for quadrotor-type aerial robots: a
  survey.
\newblock {\em Transactions of the Institute of Measurement and Control} {\bf
  2016}, {\em 38},~529--554.

\bibitem[Bouadi \em{et~al.}(2007)Bouadi, Bouchoucha, and
  Tadjine]{bouadi2007modelling}
Bouadi, H.; Bouchoucha, M.; Tadjine, M.
\newblock Modelling and stabilizing control laws design based on backstepping
  for an {UAV} type-quadrotor.
\newblock {\em IFAC Proceedings Volumes} {\bf 2007}, {\em 40},~245--250.

\bibitem[Li \em{et~al.}(2014)Li, Li, and Geng]{li2014adaptive}
Li, S.; Li, B.; Geng, Q.
\newblock Adaptive sliding mode control for quadrotor helicopters.
\newblock  Proceedings of the 33rd Chinese control conference. IEEE,  2014, pp.
  71--76.

\bibitem[Tran \em{et~al.}(2020)Tran, Santoso, Garratt, and
  Anavatti]{tran2020distributed}
Tran, V.P.; Santoso, F.; Garratt, M.; Anavatti, S.
\newblock Distributed Artificial Neural Networks-Based Adaptive Strictly
  Negative Imaginary Formation Controller for Unmanned Aerial Vehicles in
  Time-Varying Environments.
\newblock {\em IEEE Transactions on Industrial Informatics} {\bf 2020}.

\bibitem[Alexis \em{et~al.}(2012)Alexis, Nikolakopoulos, and
  Tzes]{alexis2012model}
Alexis, K.; Nikolakopoulos, G.; Tzes, A.
\newblock Model predictive quadrotor control: attitude, altitude and position
  experimental studies.
\newblock {\em IET Control Theory \& Applications} {\bf 2012}, {\em
  6},~1812--1827.

\bibitem[Lizarralde and Wen(1996)]{lizarralde1996attitude}
Lizarralde, F.; Wen, J.T.
\newblock Attitude control without angular velocity measurement: A passivity
  approach.
\newblock {\em IEEE transactions on Automatic Control} {\bf 1996}, {\em
  41},~468--472.

\bibitem[Ghallab \em{et~al.}(2018)Ghallab, Mabrok, and
  Petersen]{ghallab2018extending}
Ghallab, A.G.; Mabrok, M.A.; Petersen, I.R.
\newblock Extending Negative Imaginary Systems Theory to Nonlinear Systems.
\newblock  2018 IEEE Conference on Decision and Control (CDC). IEEE,  2018, pp.
  2348--2353.

\bibitem[Petersen and Lanzon(2010)]{petersen2010}
Petersen, I.R.; Lanzon, A.
\newblock Feedback Control of Negative Imaginary Systems.
\newblock {\em IEEE Control System Magazine} {\bf 2010}, {\em 30},~54--72.

\bibitem[Lanzon and Petersen(2008)]{lanzon2008}
Lanzon, A.; Petersen, I.R.
\newblock Stability Robustness of a Feedback Interconnection of Systems With
  Negative Imaginary Frequency Response.
\newblock {\em IEEE Transactions on Automatic Control} {\bf 2008}, {\em
  53},~1042--1046.

\bibitem[Xiong \em{et~al.}(2010)Xiong, Petersen, and Lanzon]{xiong21010jor}
Xiong, J.; Petersen, I.R.; Lanzon, A.
\newblock A negative imaginary lemma and the stability of interconnections of
  linear negative imaginary systems.
\newblock {\em IEEE Transactions on Automatic Control} {\bf 2010}, {\em
  55},~2342--2347.

\bibitem[Lewis \em{et~al.}(2003)Lewis, Dawson, and Abdallah]{lewis2003robot}
Lewis, F.L.; Dawson, D.M.; Abdallah, C.T.
\newblock {\em Robot manipulator control: theory and practice}; CRC Press,
  2003.

\bibitem[Ortega and Spong(1989)]{ortega1989adaptive}
Ortega, R.; Spong, M.W.
\newblock Adaptive motion control of rigid robots: A tutorial.
\newblock {\em Automatica} {\bf 1989}, {\em 25},~877--888.

\bibitem[Bolandi \em{et~al.}(2013)Bolandi, Rezaei, Mohsenipour, Nemati, and
  Smailzadeh]{bolandi2013attitude}
Bolandi, H.; Rezaei, M.; Mohsenipour, R.; Nemati, H.; Smailzadeh, S.M.
\newblock Attitude control of a quadrotor with optimized {PID} controller {\bf
  2013}.

\bibitem[Thanh \em{et~al.}(2018)Thanh, Phi, and Hong]{thanh2018simple}
Thanh, H.L.N.N.; Phi, N.N.; Hong, S.K.
\newblock Simple nonlinear control of quadcopter for collision avoidance based
  on geometric approach in static environment.
\newblock {\em International Journal of Advanced Robotic Systems} {\bf 2018},
  {\em 15},~1729881418767575.

\end{thebibliography}
\bibliographystyle{chicago2}
\end{document}